\documentclass[preprint,12pt,authoryear]{elsarticle}

\usepackage{hyperref} 
\usepackage{amssymb,amsmath}
\usepackage{graphicx}
\usepackage{rotating,booktabs}
\usepackage{xcolor}
\newcommand{\BA}[1]{{\color{black}#1}}

\begin{document}

\begin{frontmatter}

\title{Accelerated Benders decomposition for bi-objective mixed integer programming}

\author[label1]{Duleabom An\corref{cor1}}
\author[label2]{Chungmok Lee}
\author[label1]{Sophie N. Parragh}
\cortext[cor1]{Corresponding author. Email: dbomahn@gmail.com}

\affiliation[label1]{organization={Institute of Production and Logistics Management/JKU Business School, Johannes Kepler University Linz},
            addressline={Altenberger Strasse 69}, 
            city={Linz},
            postcode={4040}, 
            state={UpperAustria},
            country={Austria}}
            
\affiliation[label2]{organization={Department of Industrial and Management Engineering, Hankuk University of Foreign Studies},
            addressline={81 Oedae-ro, Mohyeon-eup, Cheoin-gu},
            city={Yongin-si},
            postcode={17035},
            state={Gyeonggi-do},
            country={South Korea}}

\begin{abstract}
We investigate bi-objective mixed integer programming to which Benders decomposition can be applied, a widely used technique for exploiting problem structures in real-world applications. Our study focuses on a sustainable supply chain network design problem that aims to minimise both economic and environmental objectives. To address the bi-objective nature, we employ the weighted sum method, which requires solving similar problems repeatedly with different weight values.
To tackle this, we propose a straightforward acceleration technique that reuses Benders cuts generated in previous iterations of the weighted sum method. This approach enables the Benders decomposition master problem to be initialised with pre-generated cuts. Our computational experiments demonstrate that reusing these cuts significantly reduces computing time.
\end{abstract}





\end{frontmatter}

\section{Introduction}

In the diverse fields of engineering and industry, decision-making processes often involve conflicting objectives, giving rise to multi-objective optimisation problems. Such problems require finding a set of trade-off solutions, where the improvement of one objective cannot be achieved without a deterioration in another objective.
Furthermore, real-world applications are often large-scale, and they consequently pose significant computational limitations, making them challenging to solve within reasonable computation time.
One approach that has demonstrated its effectiveness in the single-objective domain is Benders decomposition (BD) initially proposed by \citet{Benders1962}.
Based on our current understanding, only a limited number of studies have applied the classical BD framework to the bi-objective setting.  
In the domain of scalarisation-based methods, 
\BA{\citet{mardan2019accelerated} combined weighted-sum scalarisation with Benders decomposition to address a bi-objective green closed-loop supply chain network design problem. Their acceleration strategies focus on warm-starting and problem-specific enhancements, while each scalarised problem is solved independently.}
\citet{parragh2021branch} applied Benders' idea to address the bi-objective stochastic facility location problem. The researchers integrated BD into the branch-and-bound scheme to compute the lower bound set. Their problem setting required only optimality cuts, which affected only one of the two objectives considered.
\BA{\citet{dowson2022biobjective} studied bi-objective multistage stochastic linear programming problems and combined weighted-sum scalarisation with an L-shaped decomposition method. Their methodology exploits the continuous linear programming structure and solves each scalarised problem independently.}
\citet{raith2025benders} proposed a theoretical framework for applying Benders decomposition to bi-objective linear programming and integrated it with the bi-objective simplex algorithm to generate a complete set of efficient solutions.

One potential concern in employing a scalarisation-based method is the repeated need to solve similar problems with different weight values. Therefore, we propose a simple approach to reuse the generated Benders cuts, thereby accelerating the subsequently solved problems.

\section{Accelerated Benders Decomposition for Bi-objective Mixed Integer Programming}
In this section, we present a general mathematical model for bi-objective BD and introduce an idea that accelerates the BD algorithm.

\subsection{Bi-objective Benders Decomposition}\label{BOBD}
We consider a bi-objective mixed integer programming problem defined as follows:
\[
\min
\begin{cases}
f_1 y + c_1 x\\  
f_2 y + c_2 x
\end{cases}
\quad \text{s.t.} \quad Dy + Ex \leq d, \quad y \in \mathbb{Z}^n_{\geq}, \quad x \in \mathbb{R}^l_{\geq},
\tag{BP}
\]
where \( f_1, f_2 \in \mathbb{R}^n_{\geq} \) and \( c_1, c_2 \in \mathbb{R}^l_{\geq} \) are the cost vectors associated with integer variables \( y \in \mathbb{Z}^n_{\geq} \) and continuous variables \( x \in \mathbb{R}^l_{\geq} \), respectively. The matrix \( D \in \mathbb{R}^{m \times n} \) captures the constraints associated with \( y \), and \( E \in \mathbb{R}^{m \times l} \) those associated with \( x \). The right-hand side is given by \( d \in \mathbb{R}^m_{\geq} \).

To address the bi-objective structure, we apply the weighted sum method \citep{aneja1979bicriteria}, which scalarises the objectives using a convex combination based on a weight \( \lambda \in (0,1) \). 
For a given weight \( \lambda \), the objectives can be combined into
$f^\lambda y +c^\lambda x$, where $f^\lambda = \lambda f_1 + (1-\lambda) f_2, \quad c^\lambda = \lambda c_1 + (1-\lambda) c_2.$



In BD, the original problem is divided into a master problem, which typically contains integer variables, and a subproblem, usually a linear programme, involving the remaining continuous variables. The algorithm proceeds iteratively: at each iteration, the master problem provides a candidate solution for the integer variables, which is then passed to the subproblem for evaluation. In our approach, we solve the dual subproblem, whose feasibility depends on the current master problem solution, and from which Benders cuts are derived and added to the master problem accordingly to tighten the model.

In the first stage of BD, we solve the following master problem MP$^\lambda$, which accumulates feasibility and optimality cuts derived from the subproblem solutions. The objective is defined for a fixed weight \( \lambda \) from the weighted-sum method. 
\BA{
\begin{equation*}
    \begin{aligned}
        \min_{\substack{y \in \mathbb{Z}_\geq^n,  \eta \in \mathbb{R}}} \quad & f^{\lambda} y + \eta \\
        \text{s.t.} \quad & \eta \geq  (d - Dy)^\top \pi_p, \quad \forall p \in P^\lambda \quad \text{(optimality cuts)} \\
                         & 0 \geq (d - Dy)^\top r_q , \quad \forall q \in Q^\lambda \quad \text{(feasibility cuts)} 
    \end{aligned} \tag{MP$^\lambda$}
\end{equation*}
where \( \eta \) is an auxiliary variable underestimating the subproblem’s objective, and \( P^\lambda \), \( Q^\lambda \) denote sets of generated optimality and feasibility cuts, respectively.

For a fixed integer solution \(\hat{y}\), the dual subproblem DSP$^\lambda(\hat y)$ is:
\begin{equation*}
    \begin{aligned}
        \max_{\pi \in \mathbb{R}^l_\geq } \quad & (d - D \hat{y})^\top \pi \\
        \text{s.t.} \quad & E^\top \pi \leq c^\lambda.
    \end{aligned} \tag{DSP$^\lambda(\hat y)$}
\end{equation*}
}
If DSP$^\lambda(\hat y)$ is unbounded, an extreme ray $r_q$ is obtained, leading to a feasibility cut added to MP$^\lambda$ to exclude infeasible \(y\). If DSP$^\lambda(\hat y)$ is bounded and has an optimal solution, an extreme point \(\pi_p\) is obtained, resulting in an optimality cut to tighten MP$^\lambda$ and improve the lower bound on the objective.
These cuts iteratively refine the master problem until convergence.

\subsection{Accelerated Bi-objective Benders Decomposition}\label{AccBD}

Assuming that the feasible region of DSP$^\lambda(\hat y)$ is a bounded and pointed polyhedron, the following Benders optimality cuts are valid for the master problem (MP$^{\lambda'}$) associated with a given weight $\lambda'$:
$\eta \geq (d - Dy)B_p^{-1}c^{\lambda'}  \quad \forall p \in P^{\lambda'}$,
where \( B_p(\lambda',\hat y) \) is the optimal basis matrix of DSP$^{\lambda^{\prime}}$.

The feasibility cuts generated in MP$^\lambda$ for a given weight remain valid for other weights because the extreme rays of DSP$^\lambda(\hat y)$ are determined by the constraint matrix and do not depend on the cost vector. When the weight changes, only the right-hand side of the constraints in DSP$^\lambda(\hat y)$ varies, while the left-hand side remains unchanged. Therefore, the basis associated with the extreme point of the problem with a certain weight also defines an extreme point of the problem with a different weight. This allows us to express the optimality cut coefficient vector as \( \pi_p = c^\lambda B_p^{-1}(\lambda, \hat y) \). Based on this observation, after obtaining the basis, we add pre-generated Benders cuts to MP$^\lambda$ before commencing with a new weight to accelerate the algorithm speed.

\BA{
To avoid explicitly forming and inverting the optimal basis matrix \( B_p \), we employ CPLEX's internal basis representation. After solving DSP$^\lambda(\hat y)$, the optimal basis is retained internally. Using the CPLEX API, we first call \texttt{CPXXgethead} to retrieve the basis header, which identifies the order of basic variables. Then, we use \texttt{CPXXbtrans} to solve linear systems of the form \( \pi B_p = c^\lambda \), without computing \( B_p^{-1} \) directly. 
Specifically, we solve
$\pi_1 B_p = c_1 \quad \text{and} \quad \pi_2 B_p = c_2$,
obtaining two coefficient vectors \( \pi_1 \) and \( \pi_2 \). 
In practice, these vectors may contain negative components, so any cut based on their convex combination is valid only when the resulting vector is non-negative.
For the detailed procedure to compute the optimality cut coefficient vector \( \pi^\lambda \), we refer to \citet{raith2025benders}.
}

Compared to our method, 
\BA{\citet{raith2025benders} develop a theoretical framework for applying Benders decomposition to bi-objective linear programming by integrating it with a simplex algorithm, enabling the reuse of cuts within a single search process.}
In contrast, our method tackles mixed integer programming problems, where such integration with the simplex method is not directly applicable. Instead, we iteratively find optimal extreme points within a dichotomic search framework, which enables the reuse of generated cuts in subsequent searches and gives rise to a process resembling a branch-and-cut scheme.
\BA{Furthermore, the ``weighted optimality cuts'' in \citet{raith2025benders} should be generated from two ``extreme'' weight cases $\lambda=0$ and $\lambda=1$ by solving the respective Benders subproblems. On the other hand, our approach can generate Benders cuts for a new weight $\lambda$ from different weights, not necessarily either $\lambda=0$ or $\lambda=1$, because we store the bases $B_p$ that can yield a vaild optimality cut for \emph{any} given weight.}
\BA{While this dichotomic search systematically explores the objective space, it cannot guarantee finding all efficient solutions of the bi-objective mixed integer programming problems. Instead, it produces a representative subset of efficient solutions corresponding to the weights considered.}


\section{Computational Study}\label{sec:comp}
We evaluate our accelerated Benders decomposition algorithm (AB) in comparison with the classical Benders decomposition approach (CB). 
\BA{In CB, each weighted sum problem is solved independently from scratch, and no cuts, either optimality or feasibility, are reused across iterations. By contrast, AB reuses optimality cuts generated in previous weighted sum problems.}
To integrate BD (both CB and AB) into the bi-objective dichotomic search scheme, we first obtain two extreme solutions using lexicographic optimisation. The weight is then computed based on the objective values, following the method described in \citet{aneja1979bicriteria}.
\BA{For both CB and AB, optimality and feasibility cuts are generated during each iteration of the search. For AB, the optimal bases obtained from each subproblem are stored throughout the dichotomic search for use in subsequent iterations. From the second iteration onward, the previously generated optimality cuts are re-evaluated under the current weight. A cut is considered violated if its evaluated value exceeds the current master problem value. Among the violated cuts, a fixed number are selected in ascending order of violation and added to the master problem.}
This strategy was adopted after preliminary experiments comparing different approaches, including the most violated optimality cuts and randomly selected optimality cuts.
\BA{The number of pre-generated cuts added at each iteration was determined empirically through preliminary experiments. Cuts were increased in increments of 10, up to 100, and 40 cuts per iteration were selected (with a maximum of 40 cuts checked for violation), as this provided the best computational performance in terms of total solution time. Feasibility cuts were not reused across iterations, since initial experiments indicated that their reuse did not have a noticeable impact on computational performance.}

\subsection{Test instances}
We test the proposed scheme on a bi-objective mixed integer sustainable supply chain network design (SCND) problem adapted from \citet{eskandarpour2021multi}. The original setting is scaled down by reducing the number of facilities and products.

The network comprises four layers: suppliers ($I$), plants ($J$), distribution centres ($K$), and customers ($L$), defined on a directed graph $G = (V, A)$ with $V = I \cup J \cup K \cup L$. Products flow sequentially to meet deterministic demand. Plants and distribution centres (DCs) may use regular or green technologies, and transport options between facilities include conventional and green modes. The objectives are: (1) minimise total logistics costs (facility, production, and transport) and (2) minimise CO$_2$ emissions. Decisions involve facility locations, technology levels, transport modes, and product flows.

We generate ten small SCND instances, divided into two size classes. An instance is denoted as $(|I|,|J|,|K|,|L|)\_{\textit{index}}$, where $|I|$, $|J|$, $|K|$, and $|L|$ are the numbers of candidate suppliers, plants, DCs, and customers, respectively, and \textit{index} distinguishes instances of the same size. For all instances, $|I|=|J|$, $|K|=2|I|$, and $|L|=5|I|$, with two technologies and one product.

For conciseness, the complete mathematical formulations are provided in \ref{app1}.

\subsection{Experimental results}
All experiments were conducted on a machine equipped with an Intel Core i5-8250U CPU @ 1.60 GHz (8 cores), 16 GB of RAM. CPLEX is employed to solve the mathematical model. The dichotomic search algorithm is implemented in the programming language Julia.

\BA{We observe that in some weighted sum problems CB and AB yield different individual objective values, while still producing the same weighted sum objective value.} 
When these discrepancies affect the weights calculated during the dichotomic search, resulting in different search paths for CB and AB, we perform additional experiments using the weights generated by AB to enable a fair comparison. The corresponding results are reported in \BA{Table \ref{tb3} in \ref{app:b}.}

Table \ref{tb1} presents the following results: 
the number of efficient solutions obtained and the total number of iterations in the dichotomic search (\#ES/Iter), the number of added pre-generated \BA{optimality} cuts (\#Precut) in AB, the total number of feasibility and optimality cuts (\#Fcuts, \#Ocuts) generated during the dichotomic search, the improvement achieved by AB (denoted as imp, in \%), and the average total CPU time over five runs (denoted as CPUT, in seconds) required to complete the search.

Table \ref{tb1} shows that the numbers of feasibility and optimality cuts decrease significantly by adding pre-generated cuts from previous iterations. This naturally leads to a substantial reduction in total computing time. 
For instances such as (2,2,4,10)\_2 and (2,2,4,10)\_4, where no improvement in CPU time was observed, the time spent storing and retrieving cuts was a contributing factor. Table \ref{tb2} provides an example illustrating this behaviour.

\BA{
Table \ref{tb2} compares CB and AB for instance (2,2,4,10)\_2 across each iteration of the dichotomic search. In addition to the previously defined columns (\#Precut, \#Fcuts, \#Ocuts), the column (Weight) indicates the scalarisation weight used in the current iteration, while (subT and mpT) report the CPU time spent on solving the subproblem and master problem, respectively.
}
While pre-generated cuts reduce the number of cuts and CPU time, AB without them required more time due to the extra computation needed to obtain and store the optimal basis, as reflected in the column mpT.

Furthermore, we observe that the type of cut has a greater impact on performance than the number of cuts added, as shown in Table \ref{tb4} in \ref{app:b}. In the final iteration with weight 0.9071, adding a single cut drastically reduced both the number of generated cuts and the computation time.
\BA{
Recall that we simply choose, at most, a predetermined number of violated Benders cuts from the pool of bases obtained from previous dichotomic search iterations. Recent studies have demonstrated that adding ``stronger'' Benders cuts can enhance the performance of the Benders algorithm. For example, \cite{seoClosestBendersCut2022} proposed a reformulation technique of the Benders subproblem to identify ``closer'' Benders cuts to the relatively feasible solution. It is straightforward to implement a more ``advanced cut selection'' scheme with the proposed approach. Nonetheless, we would like to present this line of investigation as an interesting future research topic.
}

\begin{table}[h]
\scriptsize
\centering 
 \setlength{\tabcolsep}{1pt}
 \begin{tabular}{cccc*{10}{r}}
\hline  
     Instance 	   & \#ES/Iter & \#Precut &  \multicolumn{2}{c}{\#Fcuts}              &&&\multicolumn{2}{c}{\#Ocuts}               &&& \multicolumn{2}{c}{CPUT (sec)} \\ 
    \cline{4-5} \cline{8-9} \cline{12-13}
	                &		    &           &	  CB	&	           AB	&	 imp(\%)&&	    CB	&	           AB	&	imp(\%)	&&	 CB  	        &	AB	&	imp(\%)	 \\ \hline
    (2,2,4,10)\_1  	&	2 / 5	&	 40 	&	 279 	&	 \textbf{254} 	&	 9.0 	&&	 24458 	&	 \textbf{4007} 	&	 83.6 	&&	 53.3       	&	 \textbf{16.5} 	&	  69.0\\
    (2,2,4,10)\_2  	&	4 / 8	&	 80 	&	 562 	&	 \textbf{473} 	&	 15.8 	&&	 75681 	&	 \textbf{64531} &	 14.7 	&&	 \textbf{112.9}	&	  133.7	&	 -18.4 \\
    (2,2,4,10)\_3  	&	3 / 6	&	 77 	&	 1277	&	 \textbf{1062} 	&	 16.8 	&&	 13531 	&	 \textbf{9553} 	&	 29.4 	&&	 48.8       	&	 \textbf{43.1} 	&	 11.8 \\
    (2,2,4,10)\_4  	&	4 / 9	&	 87 	&	 650 	&	 \textbf{516} 	&	 20.6 	&&	 47130 	&	 \textbf{38113} &	 19.1 	&&	 \textbf{106.6}	&	 107.2 	&	 -0.6 \\
    (3,3,6,15)\_1  	&	3 / 6	&	 43 	&	 757 	&	 \textbf{605} 	&	 20.1 	&&	 35339 	&	 \textbf{30669} &	 13.2 	&&	 315.9       	&	 \textbf{291.5} 	&	 7.7 \\
    (3,3,6,15)\_2  	&	4 / 9	&	 52 	&	 852 	&	 \textbf{686} 	&	 19.5 	&&	 99831 	&	 \textbf{86897} &	 13.0 	&&	 433.4 	&	 \textbf{389.1} 	&	 10.2 \\
    (3,3,6,15)\_3  	&	4 / 9	&	 58 	&	 1857	&	 \textbf{1665} 	&	 10.3 	&&	 46510 	&	 \textbf{29947} &	 35.6 	&&	 308.1 	&	 \textbf{269.3} 	&	 12.6 \\
    (3,3,6,15)\_4  	&	4 / 9	&	 19 	&	 761 	&	 \textbf{636} 	&	 16.4 	&&	 303729 &	 \textbf{1522} 	&	 99.5 	&&	 921.8 	&	 \textbf{19.2} 	&	 97.9 \\
    (3,3,6,15)\_5  	&	4 / 9	&	 44 	&	 798 	&	 \textbf{749} 	&	  6.1   &&	   8676 &	 \textbf{7843} 	&	 9.6  	&&	   210.6 	&	 \textbf{199.7} 	&	  5.2\\
    (3,3,6,15)\_6  	&	5 / 10	&	 88 	&	 1455	&	 \textbf{1277} 	&	 12.2 	&&	 71051 	&	 \textbf{55802} &	 21.5 	&&	 347.9 	&	 \textbf{317.9} 	&	 8.6 \\
    \hline
\end{tabular}
\caption{Performance comparison between classical Benders (CB) and accelerated Benders (AB). Values in bold denote the better result between CB and AB for each metric.} 
\label{tb1}
\end{table}

\begin{table}[h] 
 \centering \setlength{\tabcolsep}{1pt}
 \begin{tabular}{cccrrrccrrcccccrr}
    \hline 
    Weight && \#Precut && \multicolumn{2}{c}{\#Fcuts}  &&& \multicolumn{2}{c}{\#Ocuts} &&&\multicolumn{2}{c}{subT (sec)} &&\multicolumn{2}{c}{mpT (sec)} \\ 
    \cline{4-6} \cline{9-10} \cline{12-14} \cline{16-17}
		&&		&&	CB	&	AB	&&&	CB	&	AB	&&&	CB	&	AB	&&	CB	&	AB	\\ \hline
0.9373	&&	0	&&	61	&	61	         &&&	51	&	51	        &&&	$\sim$0	&	$\sim$0	&&	       0.2	&	0.2	\\
0.8815	&&	0	&&	65	&	65	         &&&	713	&	713	        &&&	0.2	&	0.2	&&	       1.2	&	1.4	\\
0.8233	&&	0	&&	85	&	85	         &&&	861	&	861	        &&&	0.2	&	0.2	&&	       1.3	&	1.8	\\
0.6997	&&	0	&&	65	&	65	         &&&	15197&	15197	    &&&	2.6	&	2.7	&&	       22.4	&	32.0	\\
0.8855	&&	40	&&	62	&	\textbf{21}	 &&&	334	&\textbf{0}     &&&	0.1	&	\textbf{$\sim$0}&&	0.6	&	\textbf{0.1}	\\
0.9079	&&	0	&&	80	&	80	         &&&	688	&	688	        &&&	0.1	&	0.2	&&	         1.1	&	1.5	\\
0.8746	&&	40	&&	74	&	\textbf{26}	 &&&	14314&\textbf{3498} &&&	2.5	&	\textbf{0.6}&&	20.7	&	\textbf{7.2}	\\
0.9316	&&	0	&&	70	&	70	         &&&	43523&	43523	    &&&	7.3	&	7.6	&&	55.6	&	84.4	\\
    \hline   
\end{tabular}
\caption{Comparison between classical Benders (CB) and accelerated Benders (AB) on instance (2,2,4,10)\_2. Bold values indicate cases where BD’s performance improved with pre-generated cuts. 
} 
\label{tb2}
\end{table}

\section{Conclusion}
This study examined the application of Benders decomposition to bi-objective mixed integer programming problems, focusing on sustainable supply chain network design. To address the bi-objective nature, we employed the weighted-sum method, which involves solving a sequence of problems with varying weights. Recognising the computational burden of this approach, we proposed an accelerated Benders decomposition technique that reuses cuts generated in previous iterations. Computational results demonstrated that our method significantly reduces both total solution time and the number of generated cuts compared to the classical approach. Furthermore, our findings indicate that the type of reused cut has a greater impact on performance than their quantity. This suggests a promising direction for future research: investigating strategies for selecting cuts or developing adaptive mechanisms to balance preprocessing effort with overall efficiency. 
\BA{Another direction for future work is to test the method on simpler problem instances to better understand its behaviour and performance characteristics.}

\bibliographystyle{elsarticle-harv} 
\bibliography{mybib}
\section*{Acknowledgement}
This research was funded in whole, or in part, by the Austrian Science Fund (FWF) [P 31366 and P 32954]. This work was also supported by the National Research Foundation of Korea(NRF) grant funded by the Korea government (MSIT)(RS-2024-00452707). 
The authors would like to thank Markus Sinnl from the Institute of Business Analytics and Technology Transformation at Johannes Kepler University Linz for his valuable and insightful comments, which contributed meaningfully to this work.
For the purpose of open access, the author has applied a CC BY public copyright licence to any Author Accepted Manuscript version arising from this submission.

\newpage
\appendix
\section{Bi-Objective SCND Formulation}\label{app1}
This section provides a description of the sets, parameters, and decision variables used in the mathematical model of the bi-objective SCND problem that we use as a test case. 
It has originally been proposed by \citet{eskandarpour2021multi}, and based on their assumptions, we scaled down the original problem setting by reducing the number of facilities and products. 
Tables \ref{tb:datasets}-\ref{tb:dv} provide details on the sets, parameters and decision variables used in the model.
For a detailed description of each objective and constraint, we refer the reader to \citet{eskandarpour2021multi}.

\begin{table}[ht]
\small
\begin{center}
\begin{tabular}{cl} 
 \hline
 Set & Description\\ \hline 
 $A$ & Arcs connecting successive layers\\
 $I$  & Suppliers \\ 
 $J$ & Candidate plants\\ 
 $K$ & Candidate DCs\\
 $L$ & Customers\\
 $M$ & Transportation modes \\
 $T$ & Technologies\\
\hline
\end{tabular}
\caption{Sets}\label{tb:datasets}
\end{center}

\setlength{\tabcolsep}{1pt}
\centering
\begin{tabular}{cl} 
 \hline
 Parameter & Description\\ \hline 
 $c^{t}_j$ & Fixed cost of opening a facility $j \in J \cup K$ using technology $t \in T$\\ 
 $a_i$ & Unit processing cost at supplier $i \in I$\\
 $e^{t}_j$ & Unit processing cost at facility $j \in J \cup K$ with technology $t \in T$\\
 $cap_i$ & Capacity of facility $i \in I \cup J \cup K$\\
 $g^{m}_{ij}$ & Fixed cost of using transportation mode $m \in M$ along arc ($i,j$) $\in A$\\
 $v^{m}_{ij}$ & Variable transportation cost per unit on arc ($i,j$) $\in A$ using \\ & transportation mode $m \in M$ \\
 $\underbar{$V$}^{m}_{ij}$ & Minimum threshold volume for using transportation mode $m \in M$ \\ & along arc ($i,j$) $\in A$\\
  $\bar{V}^{m}_{ij}$ & Capacity of transportation mode $m \in M$ along arc ($i,j$) $\in A$\\
 $b_i$ & Unitary CO$_2$ emissions at supplier $i \in I$\\
 $q^{t}_j$ & CO$_2$ emissions per unit manufactured or warehoused at\\ & facility $j \in J \cup K$ with technology $t \in T$\\
 $r^{m}_{ij}$ & CO$_2$ emissions per unit transported along arc ($i,j$) $\in A$ \\ &using transportation mode $m \in M$\\
 $d_{l}$  & Demand from customer $l \in L$\\ 
 $J_{max}$ & Maximum number of plants that may be established\\
 $K_{max}$ & Maximum number of DCs that may be established\\
 \hline
\end{tabular}
\caption{Parameters}\label{tb:parameters}
\vspace{2em}
\begin{tabular}{ccll} 
 \hline
 Variable & Stage & Type & Description\\ \hline 
 $y^t_j$  & $1^{st}$ & $\mathbb{B}$ &  = 1 if facility $j \in J \cup K$ with technology level $t$ is open; 0 otherwise\\ 
 $u^m_{ij}$  & $1^{st}$ & $\mathbb{B}$ & = 1 if transportation mode $m \in M$ is selected for arc ($i,j$) $\in A$; \\&&& 0 otherwise \\
 $x^{m}_{ij}$ & $2^{nd}$ & $\mathbb{R}^+$ & Flow on arc ($i,j$) $\in A$ using transportation mode $m \in M$ \\
 $h^{t}_{j}$ & $2^{nd}$ & $\mathbb{R}^+$ & Amount processed with technology level $t \in T$ at facility $j \in J \cup K$\\
 \hline
\end{tabular}
\caption{Decision variables}\label{tb:dv}
\end{table}


\paragraph{Bi-objective SCND Model}~\\
The following presents the scaled-down bi-objective model for the sustainable supply chain network design problem.
To aid understanding of the subsequent dual subproblem formulation, we first present the constraints related to the first-stage master problem and express the model in standard form.

\begin{equation*}
    \begin{aligned}
        \min \sum_{j\in J\cup K}\sum_{t \in T}{c_{j}^{t}y_{j}^t}   + 
        \sum_{(i,j)\in A }\sum_{m \in M}{g_{ij}^{m}u_{ij}^m} 
        +\sum_{i\in I}\sum_{j \in J}\sum_{m\in M}{a^{m}_{i}x^{m}_{ij}} \\
        +\sum_{j\in J \cup K}\sum_{t \in T}{e^{t}_{j}h^{t}_{j}} 
        + \sum_{(i,j)\in A}\sum_{m \in M}{v^{m}_{ij}x^{m}_{ij}}
\end{aligned}  
\end{equation*}

\begin{equation*}
    \begin{aligned}
        \min \sum_{i\in I}\sum_{j \in J}\sum_{m\in M}{b^{m}_{i}x^{m}_{ij}} + \sum_{j\in J\cup K }\sum_{t \in T}{q^{t}_{j}h^{t}_{j}} + \sum_{(i,j)\in A}\sum_{m \in M}{r^{m}_{ij}x^{m}_{ij}}
    \end{aligned}
\end{equation*}

subject to:
\begin{align}
    \sum_{t\in T} y_{j}^{t} &\leq 1 && \forall j \in J \cup K \label{con1}\\
    \sum_{m\in M} u_{ij}^{m} &\leq 1 && \forall (i,j) \in A \label{con2}\\
    \sum_{j\in J} \sum_{t \in T} y_{j}^{t} &\leq J_{\max} \label{con3}\\
    \sum_{k\in K} \sum_{t \in T} y_{k}^{t} &\leq K_{\max} \label{con4}\\
    \sum_{i\in I} \sum_{m \in M} x_{ij}^{m} - \sum_{k\in K} \sum_{m \in M} x_{jk}^{m} &= 0 && \forall j \in J \label{con5}\\
    \sum_{j\in J} \sum_{m \in M} x_{jk}^{m} - \sum_{l\in L} \sum_{m \in M} x_{kl}^{m} &= 0 && \forall k \in K \label{con6}\\
    - \sum_{i\in I} \sum_{m \in M} x_{ij}^{m} + \sum_{t\in T} h_{j}^{t} &= 0 && \forall j \in J \label{con7}\\
    - \sum_{j\in J} \sum_{m \in M} x_{jk}^{m} + \sum_{t\in T} h_{k}^{t} &= 0 && \forall k \in K \label{con8}\\
    \sum_{k\in K} \sum_{m \in M} x_{kl}^{m} &\geq d_{l} && \forall l \in L \label{con9}\\
    - \sum_{j\in J} \sum_{m \in M} x_{ij}^{m} &\geq -cap_i && \forall i \in I \label{con10}\\
    - h_{j}^{t} &\geq -cap_{j} y_{j}^{t} && \forall j \in J \cup K, t \in T \label{con11}\\
    x_{ij}^{m} &\leq \bar{V}_{ij}^{m} u_{ij}^{m} && \forall (i,j) \in A, m \in M \label{con12}\\
    - x_{ij}^{m} &\leq -\underline{V}_{ij}^{m} u_{ij}^{m} && \forall (i,j) \in A, m \in M \label{con13}\\
    &y_{j}^{t} \in \{0,1\} && \forall j \in J \cup K, t \in T \\
    &u_{ij}^{m} \in \{0,1\} && \forall (i,j) \in A, m \in M \\
    &x_{ij}^{m} \geq 0 && \forall (i,j) \in A, m \in M \\
    &h_{j}^{t} \geq 0 && \forall j \in J \cup K, t \in T
\end{align}

\paragraph{Dual Subproblem of the Bi-objective SCND Model}~\\
The dual subproblem corresponding to the bi-objective SCND model is formulated as follows, where
$\lambda$ denotes the weight determined by the dichotomic search technique and $\alpha$ denotes the set of dual variables corresponding to the constraint numbers in the primal model.

\begin{equation*}
    \begin{aligned}
    \max \sum_{l \in L}d_{l}\alpha_{5l}- \sum_{i \in I}cap_i\alpha_{6i} - \sum_{j \in J\cup K}\sum_{t \in T}cap_j \hat{y}^t_j \alpha_{7j}^t\\
    + \sum_{(i,j) \in A}\sum_{m \in M}(\underbar{V}_{ij}^m
    -\bar{V}^{m}_{ij})\hat{u}^m_{ij}\alpha_{8ij}^m
    \end{aligned}
\end{equation*}
subject to:
\begin{align}
    \alpha_{1j} - \alpha_{3j} - \alpha_{6i} + \alpha_{9ij}^m &\leq \lambda(a_i + v_{ij}^m) + (1-\lambda)(b_i + r_{ij}^m) \quad \forall i \in I, j \in J, m \in M \label{con14} \\
    -\alpha_{1j} + \alpha_{2k} - \alpha_{4k} + \alpha_{10jk}^m &\leq \lambda v_{jk}^m + (1-\lambda) r_{jk}^m \quad \forall j \in J, k \in K, m \in M \label{con15} \\
    -\alpha_{2k} + \alpha_{5l} - \alpha_{11kl}^m &\leq v_{kl}^m + \lambda r_{kl}^m \quad \forall k \in K, l \in L, m \in M \label{con16} \\
    \alpha_{3j} - \alpha_{7j}^t &\leq \lambda e_j^t + (1-\lambda) q_j^t \quad \forall j \in J, t \in T \label{con17} \\
    \alpha_{4k} - \alpha_{8k}^t &\leq \lambda e_k^t + (1-\lambda) q_k^t \quad \forall k \in K, t \in T \label{con18} \\
    \alpha_{1j}, \alpha_{2k}, \alpha_{3j}, \alpha_{4k} &\text{ are free} \label{con19} \\
    \alpha_{5l} &\geq 0 \quad \forall l \in L \label{con20} \\
    \alpha_{6i} &\geq 0 \quad \forall i \in I \label{con21} \\
    \alpha_{7j}^t &\geq 0 \quad \forall j \in J, t \in T \label{con22} \\
    \alpha_{8k}^t &\geq 0 \quad \forall k \in K, t \in T \label{con23} \\
    \alpha_{9ij}^m &\geq 0 \quad \forall i \in I, j \in J, m \in M \label{con24} \\
    \alpha_{10jk}^m &\geq 0 \quad \forall j \in J, k \in K, m \in M \label{con25} \\
    \alpha_{11kl}^m &\geq 0 \quad \forall k \in K, l \in L, m \in M \label{con26}
\end{align}

\section{Additional Experimental Results}\label{app:b}
Table \ref{tb3} presents results from an additional experiment conducted on instance (2,2,4,10)\_1, where the dichotomic search scheme produced different weights. The fixed weights used are those originally derived from AB algorithm.
Table \ref{tb4} is provided to offer further insights into the experimental results. For both tables, 
\BA{columns (Weight, \#Precut, \#Fcuts, and \#Ocuts) report the scalarisation weight, number of pre-generated cuts, and the total number of feasibility and optimality cuts generated during the search. Columns (subT and mpT) indicate the CPU time spent on solving the subproblem and master problem, respectively.}
Values in bold indicate the better result between CB and AB for each metric.

\begin{table}[htp]
\setlength{\tabcolsep}{3pt}
\centering
\begin{tabular}{cccccrrccccrr}
\toprule
Weight & \#Precut & \multicolumn{2}{c}{\#Fcuts} && \multicolumn{2}{c}{\#Ocuts} && \multicolumn{2}{c}{subT (sec)} && \multicolumn{2}{c}{mpT (sec)} \\
\cmidrule{3-4} \cmidrule{6-7} \cmidrule{9-10} \cmidrule{12-13} 
& & CB & AB && CB & AB && CB & AB && CB & AB \\
\midrule
    0.8579 & 0 & 56 & 56 && 196 & 196 && 0.1 & 0.1 && 1.0 & 1.1 \\
    0.9395 & 4 & 56 & \textbf{46} && 261 & \textbf{144} && 0.1 & 0.1 && 1.1 & \textbf{0.8} \\
    0.8050 & 35 & 52 & \textbf{44} && 22327 & \textbf{3359} && 6.6 & \textbf{1.0} && 43.8 & \textbf{9.5} \\
    0.7551 & 0 & 56 & 56 && 247 & 247 && 0.1 & 0.1 && 0.9 & 1.2 \\
    0.8343 & 1 & 59 & \textbf{52} && 1427 & \textbf{61} && 0.5 & \textbf{$\sim$0} && 3.4 & \textbf{0.5} \\
\bottomrule
\end{tabular}
\caption{Comparison between classical Benders (CB) and accelerated Benders (AB) on (2,2,4,10)\_1.}
\label{tb3}

\vspace{2em}

\begin{tabular}{cccrrrrrrrrrrrrrr}
\toprule
Weight && \#Precut && \multicolumn{2}{c}{\#Fcuts} &&& \multicolumn{2}{c}{\#Ocuts} &&& \multicolumn{2}{c}{subT (sec)} && \multicolumn{2}{c}{mpT (sec)} \\ 
\cline{4-6} \cline{9-10} \cline{12-14} \cline{16-17}
		&&		&&	CB	&	AB	&&&	CB	&	AB	&&&	CB	&	AB	&&	CB	&	AB	\\ \hline
    0.9195	&&	0	&&	162	&	162	&&&	977	&	977	&&&	0.3	&	0.3	&&	6.0	&	7.2	\\
    0.9579	&&	0	&&	90	&	90	&&&	531	&	531	&&&	0.2	&	0.2	&&	3.2	&	3.9	\\
    0.9222	&&	25	&&	155	&	\textbf{145}	&&&	1661	&	\textbf{912}	&&&	0.5	&	\textbf{0.3}	&&	7.8	&	\textbf{6.1}	\\
    0.9747	&&	0	&&	165	&	165	&&&	1497	&	1497	&&&	0.5	&	0.5	&&	18.7	&	18.9	\\
    0.9671	&&	7	&&	120	&	\textbf{115}	&&&	2966	&	\textbf{1399}	&&&	0.7	&	\textbf{0.4}	&&	9.7	&	\textbf{8.5}	\\
    0.8818	&&	22	&&	162	&	\textbf{93}	&&&	\textbf{16265}	&	17709	&&&	\textbf{3.7}	&	4.0	&&	\textbf{57.5}	&	78.9	\\
    0.8534	&&	0	&&	183	&	183	&&&	2861	&	2861	&&&	0.8	&	0.9	&&	16.8	&	20.9	\\
    0.8913	&&	0	&&	134	&	134	&&&	27990	&	27990	&&&	6.0	&	6.6	&&	111.8	&	150.4	\\
    0.8717	&&	33	&&	139	&	\textbf{117}	&&&	895	&	\textbf{593}	&&&	\textbf{0.3}	&	0.4	&&	\textbf{5.4}	&	10.0	\\
    0.9071	&&	1	&&	145	&	\textbf{73}	&&&	15408	&	\textbf{597}	&&&	3.3	&	\textbf{0.2}	&&	57.0	&	\textbf{5.5}	\\
\hline   
\end{tabular}
\caption{Comparison between classical Benders (CB) and accelerated Benders (AB) on instance (3,3,6,16)\_6.}
\label{tb4}
\end{table}
\end{document}